# Numerical simulation of resistance furnaces by using distributed and lumped models

Alfredo Bermúdez [1,2*], Dolores Gómez [1,2] and David González-Peñas [1]

1* Department of Applied Mathematics, Universidade de Santiago de Compostela, Lope Gómez de Marzoa, s/n. Campus Vida, Santiago de Compostela, E-15782, Spain.

2 Galician Centre for Mathematical Research and Technology, CITMAga, Constantino Candeira s/n. Campus Vida, Santiago de Compostela, E-15782, Spain.

*Corresponding author(s). E-mail(s): david.gonzalez.penas@usc.es;
Contributing authors: mdolores.gomez@usc.es; alfredo.bermudez@usc.es;

## Abstract

This work proposes a methodology combining distributed and lumped models to simulate the current distribution in an indirect heat resistance furnace and, in particular, to compute the current to be supplied in order to obtain a desired power. The distributed model is a time-harmonic eddy current problem which has been numerically solved by using a finite element method. The lumped model is based on the computation of a reduced impedance associated to an equivalent circuit model. The effectiveness of the method has been assessed by numerical simulations and plant measurements. The good correlation between the results reveals this approximation is well-suited in order to aid the design and improve the efficiency of the furnace in a short-time.

## 1 Introduction

Numerical simulation is now a tool of paramount importance for engineering design. This is particularly true in the field of electromagnetism in which Alain Bossavit's contribution has been essential for the introduction and analysis of numerical methods for solving Maxwell's equations (see, for instance, among his numerous publications the reference book [1]. The present paper deals with a good example of engineering application of this methodology to simulate the behaviour of a resistance furnace.

In the last years, indirect resistance heating has found an extending application in different manufacturing processes, especially in metallurgical, ceramic, electronic, glass and semiconductor industries ([2], [3]). In metallurgy, one of their main advantages is to produce uniform energy distribution over the workpiece, which is essential in some industrial operations.

This paper focuses on electric indirect heat resistance furnaces devoted to metal purification. In these devices, heating essentially occurs because the current supplied to furnace passes through the resistance, where electric energy is transformed into thermal energy and then heat is transferred to the load by radiation, convection and/or conduction. For furnaces working at high temperatures, radiation is the major mode of heat transfer. Usually, alternating current is used and then heat is also produced, but to a much lesser extent, due to the Joule effect produced by the eddy currents induced in the conducting parts of the domain.

This is an intensive energy consuming industrial process, and the efforts made to optimise its operation are great. In general, resistance heating is used when demand patterns are not appropriate to induction furnaces and hence cannot be considered as a competitive technique [4]. For instance, in metal melting processes, resistance furnaces have the advantage of allowing better temperature control, which is achieved by feedback with the power supply [5].

The global design requirements generally involve different and complex models because they are subjected to several physical phenomena: electromagnetic, heat, structural, fluid dynamics...Therefore, it is interesting to address their study from a multiphysics perspective. In the last years, numerical simulation has revealed as an important tool for predicting furnaces behaviour, since it allows to perform changes in silico, thus avoiding trial-error in plant operations and, in particular, those concerning to the feedback procedures with the power supply. Usually, this approach requires the use of modelling methodologies aimed at reaching a compromise between accurate results and reasonable computational times. Even though resistance furnaces are generally designed for a particular application, it is quite common modifying an existing furnace to operate at a different temperature range depending on the needs. This can be done by adjusting the power supplied to the furnace, which, in its turn, involves estimating the input current. The need for quick models for control algorithms and real time simulation is therefore important.

The prediction of the current distribution inside the resistor is a non-trivial matter, which can be only properly performed by numerical solution of the underlying electromagnetic models. This has already been done with other kind of classic metallurgical furnaces (e.g.: induction furnaces [6], [7], arc furnaces [8], [9], or electron-beam furnaces [10], [11]). Contrary to all these technologies, the bibliography relative to the mathematical modelling of industrial examples of indirect resistance furnaces is practically non-existent, which is a handicap when evaluating a new product.

In most cases, the heater has not cylindrical symmetry and as a consequence three dimensional simulations are needed, which are long and difficult. Moreover, at the design stage, it is in general mandatory to know the response to different inputs, so a single simulation is not enough.

In this article we introduce a methodology which requires a single numerical electromagnetic simulation to predict the response of the furnace to different current inputs. More precisely, we show how, starting from this initial numerical simulation, it is possible to obtain an equivalent lumped model for precisely estimating the current to be supplied to the system in order to obtain the desired power.

The use of lumped models is something usual in industrial design (e.g., induction furnaces in [12], [13], electrochemical cells in [14] or electromagnetic actuators in [15]). In our case, the novelty lies in the type of furnace and also in the fact that the method proposed can be applied to any resistance furnace powered by a three-phase alternating current, and whatever may be the geometry of the resistance and the rest of the elements composing the device.

## 2 Statement of the problem

### 2.1 Geometry of the furnace

We consider a resistance furnace as the one sketched in Fig. 1. It consist of a stainless-steel chamber enclosing a resistive heater and a workpiece placed over the heater. The workpiece is formed by a hemispherical crucible and the load within, which is the metal to be heated. Both the crucible and the resistance are surrounded by air and isolated materials. Moreover, the chamber has an internal

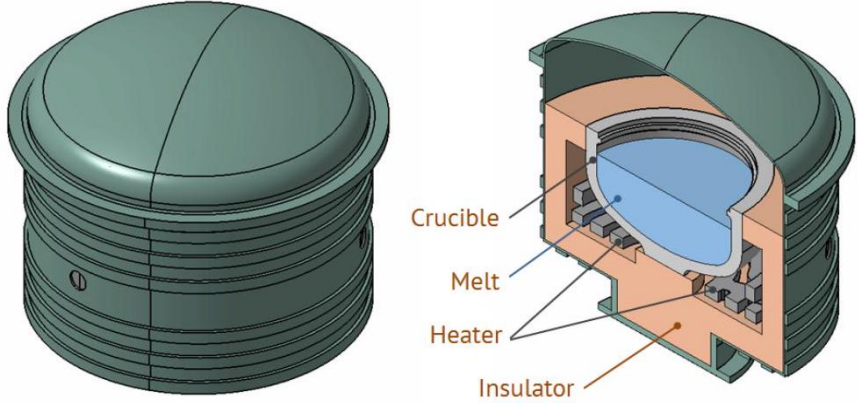

**Fig. 1**. Furnace geometry. Isometric view (left) and section (right).

insulating lining which helps to avoid heating losses.

The system is fed with a 3-phase alternating current (AC) for better efficiency (in some applications, this guarantee the good stirring of the material once it is melted).

The current is supplied to the resistor through 3 electrodes traversing the chamber up to the supply circuit. Each of them is powered by alternating current of the same amplitude but with different phases. These currents are supposed to be the known input data.

### 2.2 Mathematical model

For modelling purposes, we will introduce some notations. Let $\Omega$ be a three-dimensional bounded domain consisting of two disjoint parts, $\Omega_c$ and $\Omega_d$, corresponding to conductors and dielectrics. We will denote by $\Omega_R$ the component of $\Omega_c$ consisting of the resistive heater and its terminals. Moreover, we will denote $\Gamma$ the boundary of $\Omega$; $\bar{\Gamma}_R := \partial\Omega_R \cap \partial\Omega$ will stand for the outer boundary of $\Omega_R$ and $\bar{\Gamma}_d := \partial\Omega_d \cap \partial\Omega$ for that of the dielectric domain. Finally, **n** represents any outward unit normal vector to a given surface. The computational domain is sketched in Fig. 2 according to the previous notations. The goal of this study is to compute the distribution of current density inside the heater and the Joule effect.

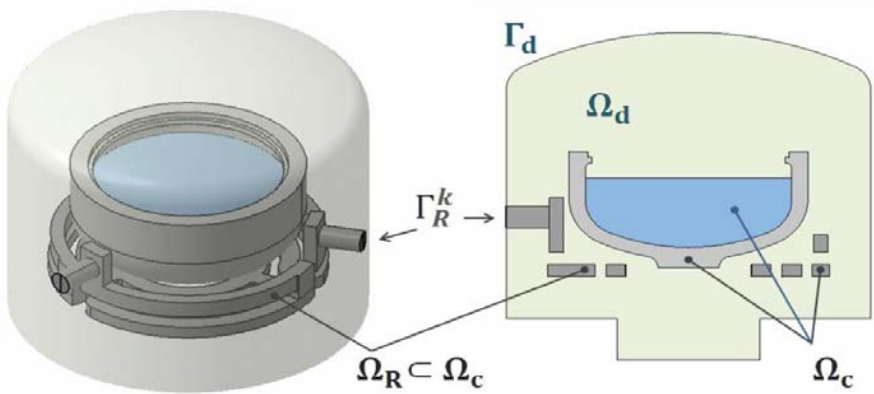

**Fig. 2.** Computational domain and notations. Isometric view (left) and section (right).

The electromagnetic model is based on the well-known time-harmonic eddy current model and it will be numerically solved by using a finite element method. Since the current source is alternating and all materials are supposed to have a linear electromagnetic behaviour, we can use a time-harmonic approach. Thus, all of the fields involved in the Maxwell system have the form

$$\mathcal{F}(\mathbf{x}, t) = \mathrm{Re}\left[\mathbf{F}(\mathbf{x})e^{i\omega t}\right] \tag{1}$$

where $t$ is time, $\mathbf{x} \in \mathbb{R}^3$ is the space position, i is the imaginary unit, $\mathbf{F}(\mathbf{x})$ is the complex amplitude (or phasor) of the field $\mathcal{F}$ and $\omega$ the angular frequency, $\omega = 2\pi f$, $f$ being the frequency of the alternating current. For industrial cases with low and moderate frequencies, quasi-static assumption applies and the term corresponding to the displacement current in the Ampere's law can be neglected. Moreover, taking into account that the electric field is not required in non-conducting materials, the time-harmonic eddy current model leads to solve the following equations defined in $\Omega$ :

$$\mathrm{curl}\ \mathbf{H} = \mathbf{J}, \tag{2}$$

$$i\omega\mathbf{B} + \mathrm{curl}\ \mathbf{E} = \mathbf{0}, \tag{3}$$

$$\mathrm{div}\ \mathbf{B} = 0, \tag{4}$$

where $\mathbf{H}$, $\mathbf{J}$ and $\mathbf{E}$ are the complex amplitudes associated with the magnetic field, the current density and the electric field, respectively. For more details about these equations see, for instance, the book of [16] or the book of [17]. To obtain a closed system, we add the constitutive law $\mathbf{B} = \mu\mathbf{H}$, $\mu$ being the magnetic permeability, and the Ohm's law, $\mathbf{J} = \sigma\mathbf{E}$, $\sigma$ being the electric conductivity which is greater than zero in conductors and null in dielectrics.

The model must be completed with suitable boundary conditions, and we consider the following ones:

$$\mu\mathbf{H} \cdot \mathbf{n} = 0 \text{ on } \Gamma \tag{5}$$

$$\mathbf{E} \times \mathbf{n} = \mathbf{0} \text{ on } \Gamma_R \tag{6}$$

Condition (5) implies that the magnetic field is tangential to the boundary of the chamber, whereas (6) means that electric current enters the domain perpendicular to the cross section of the electrodes. Notice that $\Gamma_R$ consists of the three disjoint connected components $\Gamma_R^1$, $\Gamma_R^2$ and $\Gamma_R^3$ corresponding to the top of the electrodes where the current (or the voltage) will be prescribed. In particular, we will prescribe the current in two of the electrodes, i.e,

$$\int_{\Gamma_R^k} \mathbf{J} \cdot \mathbf{n} \mathrm{d}S = -I_k \ k = 1,2, \tag{7}$$

and on the third one we will impose a null potential representing the ground

$$V_3 = 0 \text{ on } \Gamma_R^3 \tag{8}$$

The complex functions $I_k = I_k e^{i\iota_k}$ in (7) are the phasors corresponding to the harmonic time-dependent signals $\mathcal{I}_k(t) = I_k \cos(\omega t + \iota_k)$, $k = 1,2,3$, which are the real measurements at the terminals. Similarly, the complex functions $V_k = V_k e^{i\epsilon_k}$ are the phasors corresponding to the signals $\mathcal{V}_k(t) = V_k \cos(\omega t + \epsilon_k)$, $k = 1,2,3$.

The previous model can be approximated by using different formulations [18]. In this work we have considered the one based on magnetic vector potential/scalar electric potential, namely, the well-known A/V formulation. We recall here that the magnetic vector potential comes from equation (4) which implies that there exists a vector field $\mathbf{A}$, such that $\mathbf{B} = \mathrm{curl}\ \mathbf{A}$. On the other hand, the electric scalar potential V in the conducting domain comes from Faraday's law curl $(\mathbf{E} + i\omega\mathbf{A}) = \mathbf{0}$ which implies that $\mathbf{E} + i\omega\mathbf{A} = -\mathrm{grad}\ V$. The gauge condition div $\mathbf{A} = 0$ joint with suitable boundary conditions are added to ensure the uniqueness of the magnetic vector potential.

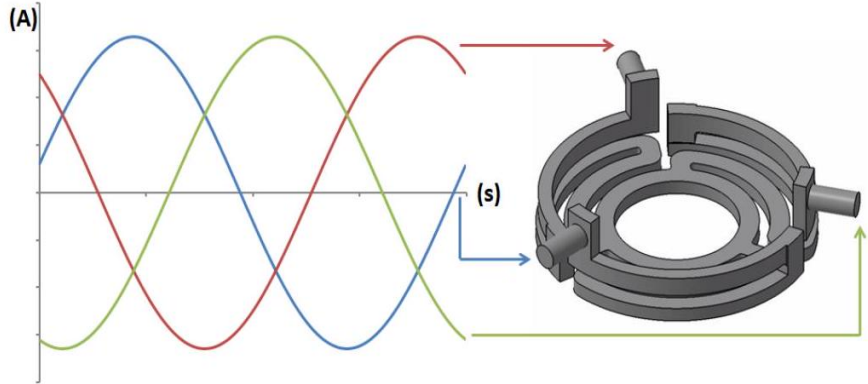

**Fig. 3.** Example of supplied intensity $(\mathcal{I}_k(t), k = 1,2,3)$.

Summarizing, the problem reads as follows:

Given complex numbers $I_k, k = 1,2$, find a vector field $\mathbf{A}$ defined in $\Omega$, and a scalar field V defined in $\Omega_c$ and constant in $\Gamma_R^1, \Gamma_R^2$, such that

$$\sigma(i\omega\mathbf{A} + \mathrm{grad}\ V) + \mathrm{curl}\left(\frac{1}{\mu}\mathrm{curl}\ \mathbf{A}\right) = \mathbf{0} \text{ in } \Omega,$$

$$\mathrm{div}\ \mathbf{A} = 0 \text{ in } \Omega,$$

$$\mathbf{A} \times \mathbf{n} = 0 \text{ on } \Gamma,$$

$$\sigma(i\omega\mathbf{A} + \mathrm{grad}\ V) \cdot \mathbf{n} = 0 \text{ on } \partial\Omega_c \setminus \Gamma_R,$$

$$V = 0 \text{ on } \Gamma_R^3,$$

$$\int_{\Gamma_R^k} \sigma(i\omega\mathbf{A} + \mathrm{grad}\ V) \cdot \mathbf{n}\mathrm{d}S = -I_k, k = 1,2.$$

### 3 Equivalent lumped model

As we have already mentioned, the furnace is connected to an external electric circuit which provides the input currents to the electrodes. We assume that the heater has no external contacts other than these three terminals. Moreover, we also assume that it is perfectly isolated inside the chamber and that there is no risk of leakage currents or earth faults. Under these assumptions, the circuit sketched in Figure 4 applies.

The geometry of the resistance as well as its location near the workpiece make appear different induction phenomena. As a consequence, to state an equivalent circuit for the whole furnace is, in general, a difficult task. Nevertheless, if we combine all the passive elements of the circuit within the same block and forget about its particular topology, it is possible to replace it by a multi-terminal connected network as the one sketched in Figure 5, with six associated variables: $\mathcal{I}_k, \mathcal{V}_k; k = 1,2,3$. $\mathcal{I}_k$ is the current that enters through k-th terminal and $\mathcal{V}_k$ is the potential of this terminal with respect to any potential reference. The terminal currents satisfy Kirchhoff's current law which state that the algebraic sum of the currents at any node of the network is zero, i.e.,

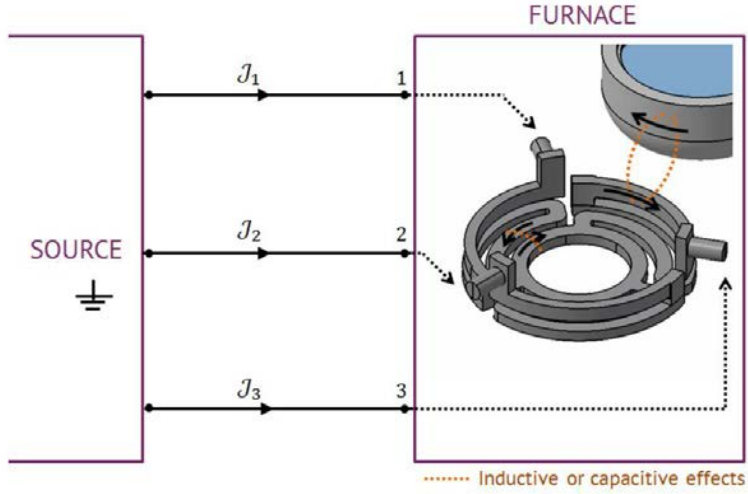

**Fig. 4.** Electrical three-phase supply via the furnace's terminals.

$$\sum_{k=1}^{3} \mathcal{I}_k = 0 \tag{9}$$

According to the general theory of electrical multi-terminal networks, which is extensively described in the books of [19] and [20], and taking into account that the network has only passive elements, currents can be expressed as a linear combination of the terminal voltages by means of the so-called indefinite admittance matrix. However, and as a consequence of (9), this is a singular matrix and, as it is shown in the articles of [21] or [22], it is not possible to obtain an impedance matrix to express the absolute voltages as a function of the currents.

Moreover, from Tellegen's Theorem [23], we deduce that the total power absorbed by the furnace network is the sum of the products $\mathcal{V}_k \mathcal{I}_k$ at the three terminals. In particular, for alternating current [24], it holds

$$S = \frac{1}{2} \sum_{k=1}^{3} V_k \bar{I}_k \tag{10}$$

where S is the total complex power absorbed by the furnace and $\bar{I}_k$ represents the conjugate of $I_k$. The active power $P$ (watts) is obtained from the real part of S. We notice that $P$ coincides with the active power $P^h$ dissipated in the heater that can be computed from the numerical simulation by using the usual formula

$$P^h = \int_\Omega \frac{|\mathbf{J(x)}|^2}{2\sigma} dV \tag{11}$$

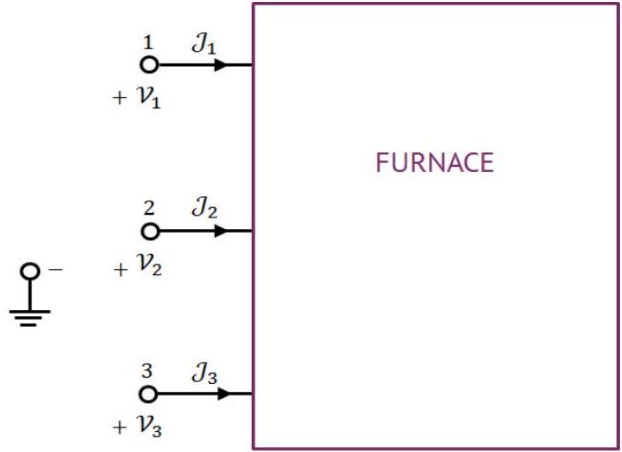

**Fig. 5.** Equivalent circuit.

Thus, we can write

$$P^{\mathrm{h}} = \mathrm{Re}\ (S). \tag{12}$$

Then, for a given current data defined by two terminal currents $I_1$ and $I_2$ we obtain the third one applying (9), i.e, $I_3 = -I_2 - I_1$ and then complex power results

$$\begin{aligned} S &= \frac{1}{2}\sum_{k=1}^{3}\ V_k\bar{I}_k \\ &= \frac{1}{2}\big(\,(\,V_1 - V_3)\bar{I}_1 + (V_2 - V_3)\bar{I}_2\big) \end{aligned}$$

Let us assume that the furnace operates with a balanced three phase supply as the one shown in Figure 3, with 120° phase shift between the electrodes, $I_2 = I_1 e^{\frac{2\pi i}{3}}$. Then, the complex power results

$$S = \frac{1}{2}\big(\,(\,V_1 - V_3) + (V_2 - V_3)e^{-i(2\pi/3)}\big)\bar{I}_1. \tag{13}$$

Equation (13) expresses the current supplied to one of the terminals as a function of the voltage drops between the terminals and the total power. We now define a reduced current and a reduced voltage drop as:

$$I^{\mathrm{R}} = I_1 \tag{14}$$

$$V^{\mathrm{R}} = (V_1 - V_3) + (V_2 - V_3)e^{-i(2\pi/3)} \tag{15}$$

As this network is linear and passive, the relation between $I^{\mathrm{R}}$ and $V^{\mathrm{R}}$ is constant. So we can define a reduced impedance as

$$Z^{\mathrm{R}} = \frac{V^{\mathrm{R}}}{I^{\mathrm{R}}} \tag{16}$$

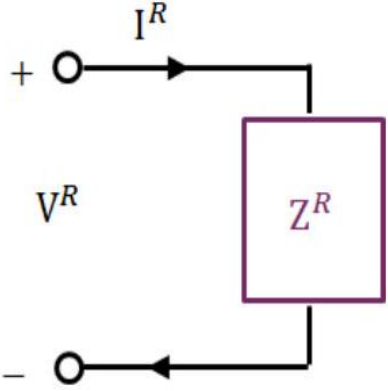

**Fig. 6.** Reduced equivalent circuit.

and use the reduced model of Figure 6 to study the furnace in terms of energy. In fact, the complex power corresponding to impedance $Z^R$ is the same as that obtained for the furnace in (13). Let us assume for simplicity that $\left|I^R\right| = I_1 = I_2 = I_3 = I$. Then, by replacing (16) in (13), we have

$$S = \frac{1}{2}Z^R I^2 \tag{17}$$

so, once we know the impedance value associated with the furnace, it is enough to know the amplitude of the terminal currents $I$ to estimate the heat dissipated in the heater:

$$P^h = \frac{1}{2}\operatorname{Re}\left(Z^R I^2\right) \tag{18}$$

Let us notice that the value of the reduced impedance can be computed by performing a single numerical simulation.

## 4 Numerical results and discussion

The numerical results presented in this section have been obtained by using a real industrial furnace the geometrical data and material properties of which are not specified here by confidentially reasons.

In the sequel, let us call 'standard' the operating conditions defined by $\left|I_1\right| = I$, $I_2 = I_1 e^{i(2\pi/3)}$ and $I_3 = I_1 e^{-i(2\pi/3)}$. The same operation conditions were applied to the furnace working in the plant in order to assess the numerical simulation and calibrate all the parameters involved. From the numerical results, it is possible to compute and represent the Joule effect in the heater (see Fig. 7) or the current density field (see Figs. 8 and 9). As a post-processing result, we can also compute the output voltage at each terminal. The results are summarized in Table 1. Finally, the lumped model of the furnace is obtained replacing the terminal voltages in (15). Thus, we get

$$V^R = V_1 + V_2 e^{-(2\pi/3)i} = 64.3238 e^{-1.4172i}(V)$$

and then

$$Z^R = \frac{V^R}{I^R} = 0.2063 e^{-3.0975i}$$

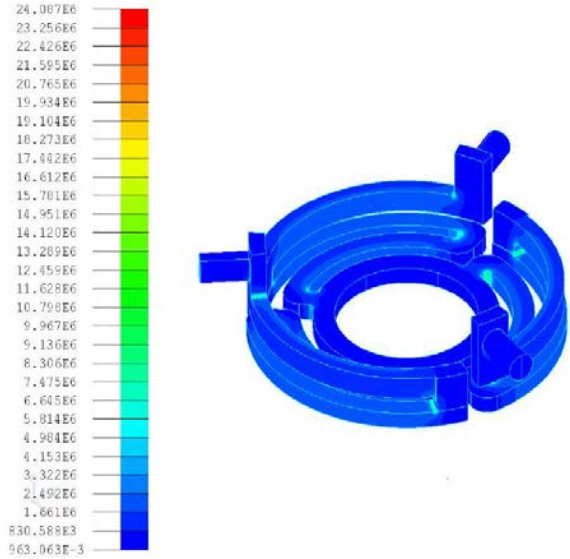

**Fig. 7.** Joule losses on the furnace heater under standard conditions (W/m³).

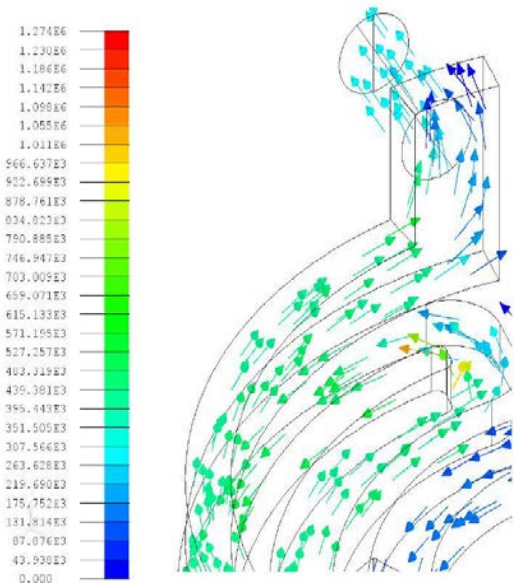

**Fig. 8.** Current density vector in the neighborhood of a terminal.

Now, the thermal dissipation in the heater may be modified in order to operate in conditions other than the standard ones with a simple procedure. Indeed, once the value $Z^R$ has been obtained, it is immediate to compute the current input which is necessary to supply to obtain the target power, by applying (18). Thus, a characteristic curve describing the dissipated power in terms of the current amplitude like the one represented in Fig. 10 can be plotted. The accuracy of this curve was validated by the company during the furnace operation. Fig. 11 represents the input currents

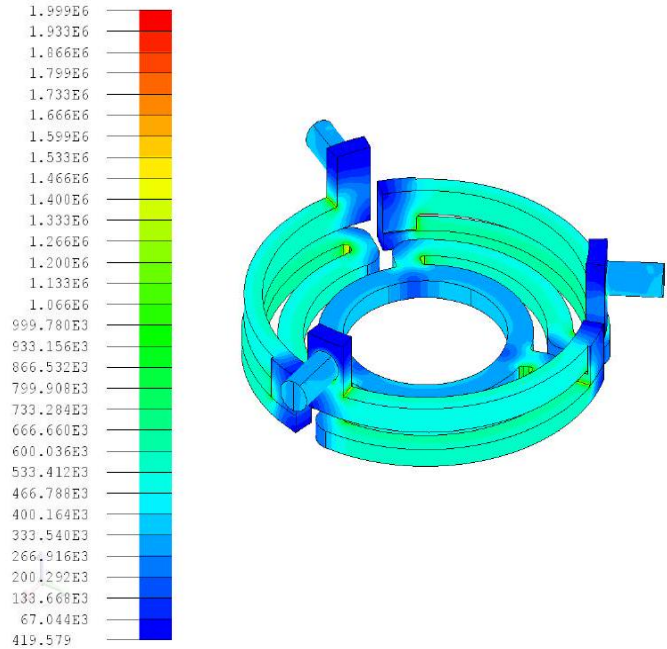

**Fig. 9.** Modulus of the current density (A/m$^2$).

**Table 1**

Terminal voltages (V, rad) and dissipated power (kW) for the standard operating conditions of current supplied (A)

| $|I^R|$ | $(V_1, \epsilon_1)$ | $(V_2, \epsilon_2)$ | $(V_3, \epsilon_3)$ | $P^h(\Omega)$ |
|---|---|---|---|---|
| 3116.79 | (37.14, −0.89) | (37.13, 0.15) | (0,0) | 100.15 |

**Table 2**

Terminal voltages (V, rad) and dissipated power (kW) as a function of current supplied (A)

| $|I^R|$ | 3266.55 | | 3411.80 | | 3939.61 | |
|---|---|---|---|---|---|---|
| $(V_1, \epsilon_1)$ | (38.92, | −0.89) | (40.65, | −0.89) | (46.94, | −0.89) |
| $(V_2, \epsilon_2)$ | (38.92, | 0.15) | (40.66, | 0.15) | (46.94, | 0.15) |
| $(V_3, \epsilon_3)$ | | (0,0) | | (0,0) | | (0,0) |
| $P^h(\Omega)$ | 109.98 | | 119.98 | | 159.97 | |

corresponding to other desired power dissipations (110, 120 and 160 kW), obtained directly from the reduced lumped model. Just to validate these values, in Table 2 we have summarized the currents and voltages obtained from the corresponding numerical simulations.

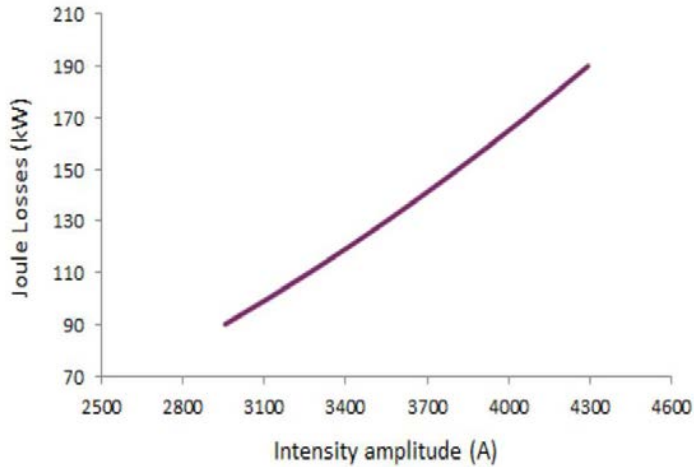

**Fig. 10.** Characteristic curve for the furnace.

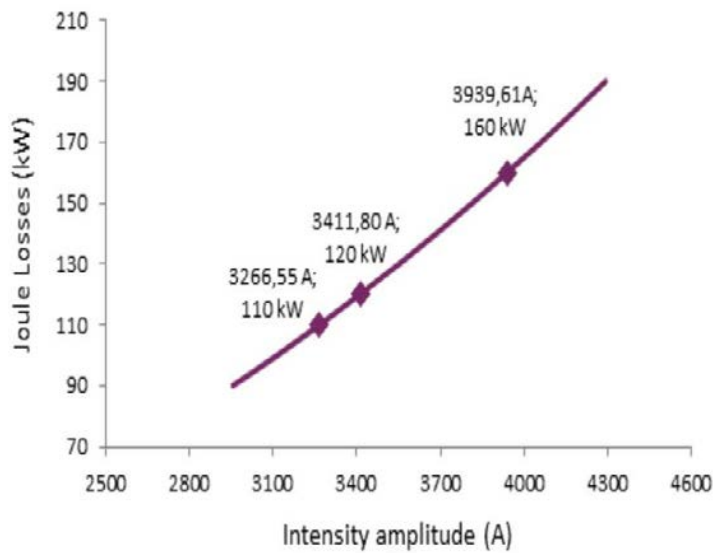

**Fig. 11**. Selection of operating points using the characteristic curve.

# 5 Some remarks

### 5.1 Reduced impedance method versus heater equivalent circuit

As an alternative to the reduced impedance procedure explained in the previous section, one could think of replacing the furnace resistance by an equivalent electrical circuit (see Fig. 12). This presents two main drawbacks. The first one is that the geometry of the resistance should be simple enough to be able to do that calculation by hand. Notice that this would be not necessary if the reduced impedance method is applied. In fact, for applying this method one should only make measurements at the terminals of the furnace in order to know its voltages. This would be not possible at the design stage but when the furnace is already operating in the plant it is very easy. On the other hand, this equivalent electrical circuit would not take into account possible induction phenomena between the furnace pieces. Although in the example shown in the present paper this phenomenon has little influence on the calculation of the impedance (it is almost 100% resistive), this may not be the case when considering other kind of furnaces operating at higher frequencies, as for instance, induction furnaces. However, the methodology based on the reduced impedance could also be easily adapted to this case: in a first step an initial numerical simulation must be done and then, in a second step and from the results of the numerical simulation, the reduced impedance can be computed.

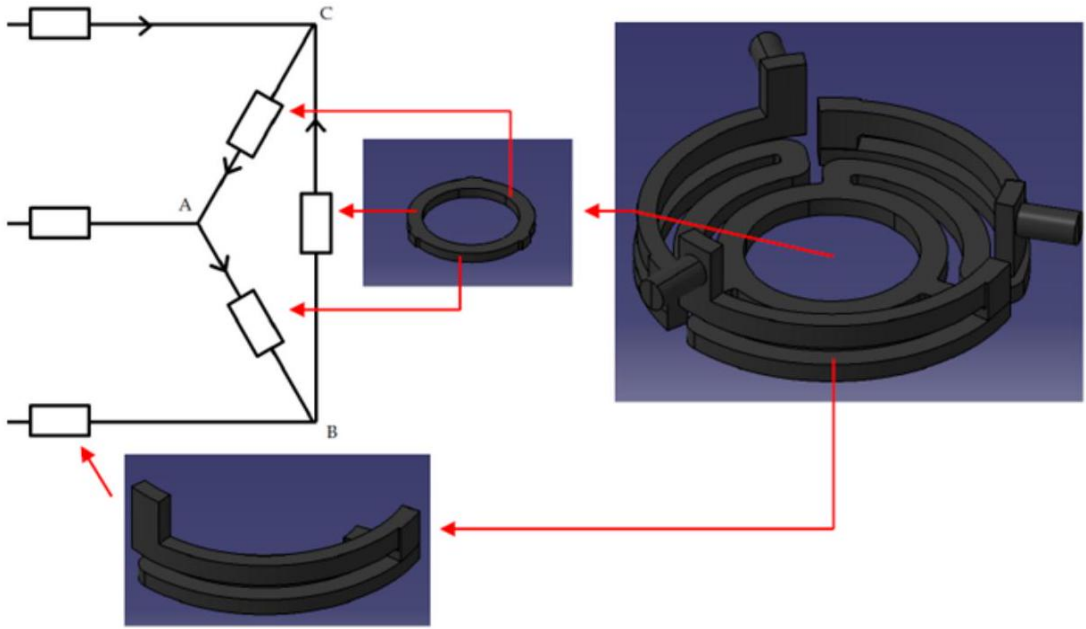

**Fig. 12.** Heater equivalent circuit construction

### 5.2 Multi-terminal network versus quadripole

Let us remember that the electric current is supplied to the furnace through 3 electrodes that go from the resistance to the power supply traversing the furnace housing. In general, this current is the only data and neither the internal geometry of the furnace nor the "upstream" circuit associated to the current supply are known. As explained in previous sections, a voltage with respect to a random reference is associated to each of the 3 terminals, so that the equivalent circuit is a multi-terminal network like the one shown in Fig. 5. Remember also that the reduced impedance method proposed in this paper is derived starting from this multi-terminal network. For this kind of circuits it is possible to write an intensity-voltage relationship in terms of the admittance matrix so that $[I] = [Y][V]$. Since $[Y]$ is a singular matrix, it is not possible to express the absolute potentials of each terminal as a function of the currents. Nevertheless, it would be possible to compute an impedance matrix for this multi-terminal network, but previously transforming it into a quadripole (see Fig. 13).

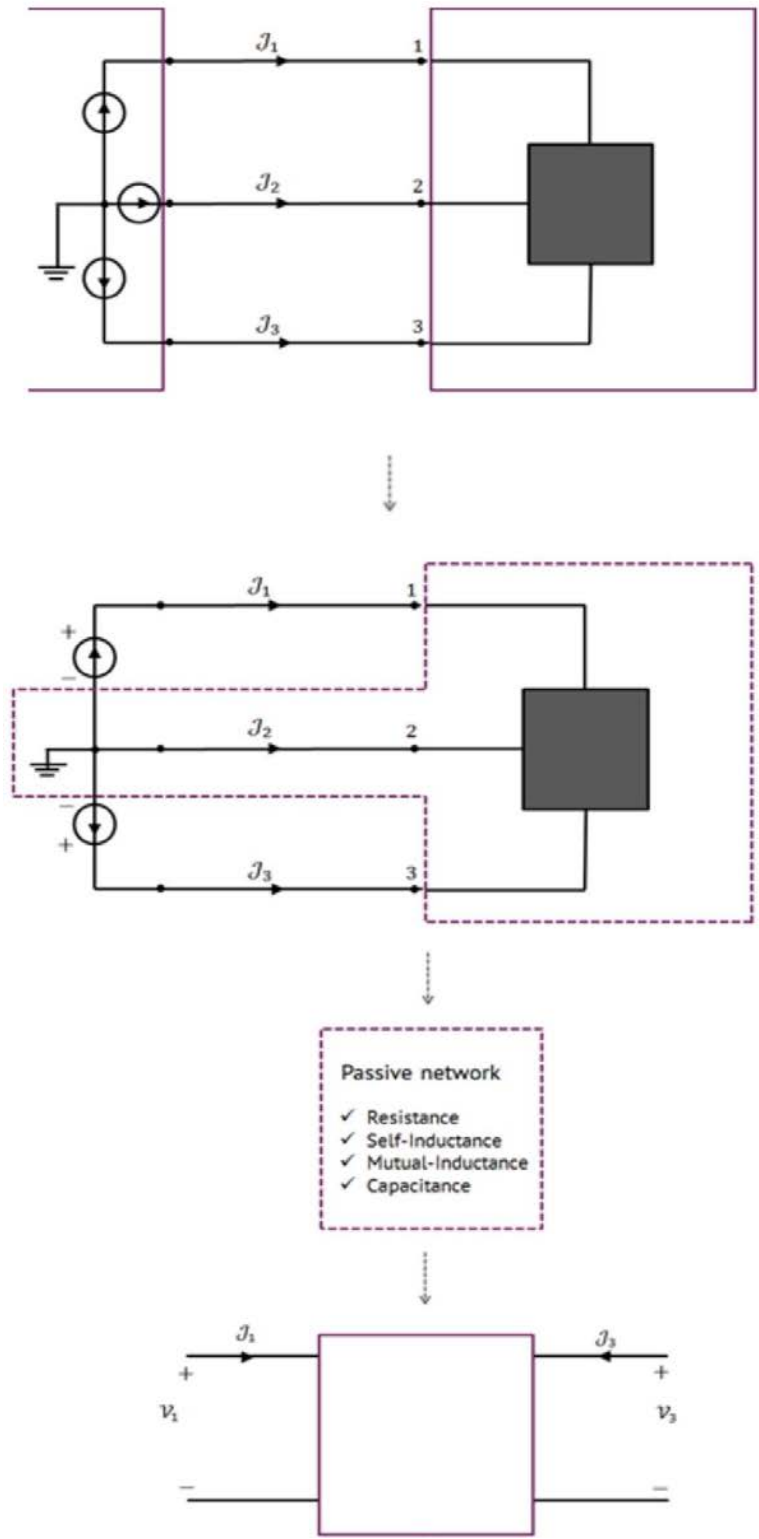

**Fig. 13.** Transforming a multi-terminal network into a quadripole

This involves to carry out a series of simplifications and hypotheses about the source circuit located "upstream", which entails a certain loss of generality. In particular, the impedance and/or admittance matrices only make sense if any non-linearity is discarded. This restriction does not appear in the reduced impedance method explained above, since only the Tellegen's theorem and Kirchhoff's laws are used and that is why the quadripole option is not considered here.

## 6 Conclusion

This paper presents a method combining a distributed and a lumped model to determine the electromagnetic behaviour of a three-phase indirect resistance furnace. The method reveals to be particularly useful at the design stage of the installation. From an initial single numerical simulation (which takes into account induction phenomena), a reduced impedance associated to the furnace is computed. From this value, it is possible to construct a curve relating the intensity supplied with the thermal power dissipated by the furnace, without carrying out more numerical simulations and avoiding long and costly trial-error procedures in plant. Moreover, the method is general in the sense that it is independent both of the upstream circuit and of the topology and materials that may exist inside the furnace. Some numerical results for a real industrial furnace are shown that assess the performance of the proposed methodology.

The authors acknowledge and appreciate the impact that Alain Bossavit's work has had on their research related to the mathematical analysis of numerical methods to solve electromagnetism problems and to the application to the mathematical modeling of various industrial problems.

## Declarations

- Funding: This work has been partially supported by FEDER and Xunta de Galicia (Spain) under grant GRC GI-1563 ED431C 2021/15.

- Conflict of interest: The authors have no conflicts of interest to declare.

- Authors' contributions: The authors contributed equally to this work.